\documentclass[12pt,leqno]{article}

\usepackage{palatino}
\usepackage{amsmath}
\usepackage{amssymb}
\usepackage{amsthm}
\usepackage{mathrsfs}
\usepackage{epsfig}
\usepackage{tabls}
\usepackage[dvips]{color}

\newtheorem{theorem}[equation]{Theorem}
\newtheorem{lemma}[equation]{Lemma}
\newtheorem{corollary}[equation]{Corollary}
\newtheorem*{sortingalg}{The $\om$-Sorting Algorithm}
\newtheorem*{exchange}{The Exchange Property}
\newtheorem*{lifting}{The Lifting Property}
\newtheorem*{gluing}{The Gluing Property}

\theoremstyle{definition}
\newtheorem{definition}[equation]{Definition}
\newtheorem*{warning}{Warning}
\newtheorem*{weakorder}{Weak Order}
\newtheorem*{bruhatorder}{Bruhat Order}
\newtheorem{example}[equation]{Example}

\newcommand{\eq}{\refstepcounter{equation}{\rm (\theequation)}}

\numberwithin{equation}{section}

\def\R{{\mathbb R}}
\def\Z{{\mathbb Z}}

\def\F{{\mathscr F}}
\def\C{{\mathscr C}}
\def\om{\omega}
\def\lex{{\sf lex}}
\def\omsort{{\sf sort}_\om}
\def\RW{{\sf R}}
\def\B{{\sf B}}
\def\rk{{\sf rk}}
\def\EL{{\rm EL}}
\newcommand{\abs}[1]{ \left|#1 \right| }

\title{The Sorting Order on a Coxeter Group}
\author{Drew Armstrong\footnote{Work supported by NSF grant DMS-0603567}\\
              \texttt{armstron@math.umn.edu}}

\begin{document}

\maketitle

\begin{abstract}
Let $(W,S)$ be an arbitrary Coxeter system. For each word $\om$ in the generators we define a partial order---called the {\sf $\om$-sorting order}---on the set of group elements $W_\om\subseteq W$ that occur as subwords of $\om$. We show that the $\om$-sorting order is a supersolvable join-distributive lattice and that it is strictly between the weak and Bruhat orders on the group. Moreover, the $\om$-sorting order is a ``maximal lattice'' in the sense that the addition of any collection of Bruhat covers results in a nonlattice.

Along the way we define a class of structures called {\sf supersolvable antimatroids} and we show that these are equivalent to the class of supersolvable join-distributive lattices.
\end{abstract}

\section{Introduction}
In this paper we will describe a very general phenomenon regarding reduced words in Coxeter groups. Let $(W,S)$ be an arbitrary Coxeter system and let
\begin{equation*}
S^*:=S^0\cup S^1\cup S^2\cup\cdots\cup S^\infty
\end{equation*}
denote the collection of finite and semi-infinite words in the generators. Given an arbitrary word $\om\in S^*$---called the {\sf sorting word}---let $W_\om\subseteq W$ denote the set of group elements that occur as subwords of $\om$. We identify every subword $\alpha\subseteq\om$ with the index set $I(\alpha)\subseteq I(\om)=\{1,2,3,\ldots\}$ describing the positions of its letters.

For each element $u\in W_\om$ let $\omsort(u)$ denote the reduced word for $u$ that is lexicographically first among subwords of $\om$---we call this the {\sf $\om$-sorted word} of $u$. In this way, $\om$ induces a canonical reduced word for each element of $W_\om$. We define the {\sf $\om$-sorting order} on $W_\om$ as the inclusion order on index sets of sorted words. It turns out that there always exists a {\em reduced} subword $\om'\subseteq\om$ such that the $\om'$- and $\om$-sorting orders coincide. Furthermore, the sorting order is not affected by the exchange of adjacent commuting generators in the sorting word. In the case of a finite Coxeter group $W$, we may summarize some of our results as follows:

\begin{quote}
For each commutation class of reduced words for the longest element $w_\circ\in W$ we obtain a supersolvable join-distributive lattice on the elements of the group. This lattice is graded by the usual Coxeter length $\ell:W\to\Z$ and it is strictly between the weak and Bruhat orders. Furthermore, the poset is a ``maximal lattice'' in the sense that the addition of any collection of Bruhat covers results in a nonlattice.
\end{quote}

More generally, the collection of $\om$-sorted subwords of $\om$ has a remarkable structure, related to the study of abstract convexity. Given a ground set $E$ and a collection of finite subsets $\F\subseteq 2^E$---called {\sf feasible} sets---we say that the pair $(E,\F)$ is an {\sf antimatroid} if it satisfies:
\begin{itemize}
\item[] The empty set $\emptyset$ is in $\F$,\quad and
\item[] Given $A,B\in\F$, $B\not\subseteq A$, there exists $x\in B\setminus A$ such that $A\cup\{x\}\in\F$.
\end{itemize}
Edelman showed that a lattice $P$ is {\sf join-distributive} (see Section \ref{sec:JD}) if and only if it arises as the lattice of feasible sets of an antimatroid.

We will prove an extension of Edelman's theorem based on the following concept. Let $(E,\leq_E)$ be a totally ordered ground set and let $\F\subseteq 2^E$ be a collection of feasible finite subsets. We say that $(E,\F,\leq_E)$ is a {\sf supersolvable antimatroid} if it satisfies:
\begin{itemize}
\item[] The empty set $\emptyset$ is in $\F$,\quad and
\item[] Given $A,B\in\F$, $B\not\subseteq A$, let $x=\min_{\leq_E}(B\setminus A)$. Then $A\cup\{x\}\in\F$.
\end{itemize}
It is clear that a supersolvable antimatroid is, in particular, an antimatroid. We prove that a lattice is join-distributive and {\sf supersolvable} if and only if it arises as the lattice of feasible sets of a supersolvable antimatroid. Finally, our main result states that
\begin{quote}
The collection of $\om$-sorted subwords of a given sorting word $\om\in S^*$ is a supersolvable antimatroid.
\end{quote}

The paper is organized as follows.

In Section \ref{sec:convexity} we review the concepts of antimatroid, convex geometry, and join-distributive lattice. After this we define supersolvable antimatroids and prove that they are equivalent to supersolvable join-distributive lattices.

Section \ref{sec:sorting} contains the definitions of $\om$-sorted words and the $\om$-sorting order. We give an algorithmic characterization of $\om$-sorted words and discuss how this algorithm is a generalization of classical sorting algorithms.

We prove our main results in Section \ref{sec:properties}. Namely, we show that the $\om$-sorting order is strictly between the weak and Bruhat orders; we prove that the collection of $\om$-sorted words forms a supersolvable antimatroid and hence that the $\om$-sorting order is a supersolvable join-distributive lattice; and we show that the $\om$-sorting order is constant on commutation classes of reduced words. Finally, we prove that the sorting orders are ``maximal lattices'' and discuss how this is related to the weak and Bruhat orders.

After this we discuss two important special cases. In Section \ref{sec:infinite} we consider the case when infinitely many group elements occur as subwords of the sorting word $\om$. In this case most of our results still hold. In particular, the $\om$-sorting order is still a lattice in which every interval is join-distributive and supersolvable. This is remarkable because the weak order on an infinite Coxeter group is {\em not} a lattice.

Finally, in Section \ref{sec:cyclic} we discuss our motivation for the current paper. This is the work of Reading on ``Coxeter-sortable elements.''

\section{Abstract Convexity}
\label{sec:convexity}

In this section we review the theory of ``abstract convexity'' by discussing three equivalent structures: antimatroids, convex geometries and join-distributive lattices. All of the structures here are finite but in Section \ref{sec:infinite} we will relax this condition. Except where stated otherwise, this material can be found in {\em Greedoids} \cite{KLS} by Korte, Lov\'asz and Schrader.

At the end of the section we will add the criterion of ``supersolvability'' and prove a characterization of supersolvable join-distributive lattices.

\subsection{Equivalent Structures}

\subsubsection{Antimatroids}
A {\sf set system} is a pair $(E,\F)$, where $E$ is a finite {\sf ground set} and $\F\subseteq 2^E$ is a collection of subsets---called {\sf feasible} sets. The system is {\sf accessible} if it satisfies:
\begin{itemize}
\item[] For each nonempty $A\in\F$, there exists $x\in A$ such that $A\setminus\{x\}\in\F$.
\end{itemize}
More specifically, an accessible set system is called an {\sf antimatroid} if it satisfies any of the following equivalent conditions.
\begin{lemma}{\bf \cite[Lemma III.1.2]{KLS}}
\label{lem:antimatroid}
Given an accessible set system $(E,\F)$, the following statements are equivalent:
\begin{itemize}
\item[\eq\label{eq:antimatroid}]
For all feasible sets $A,B\in\F$ with $B\not\subseteq A$, there exists $x\in B\setminus A$ such that $A\cup\{x\}\in\F$.
\item[\eq\label{eq:unions}] $\F$ is closed under taking unions,
\item[\eq] Given $A$, $A\cup\{x\}$ and $A\cup\{y\}$ in $\F$, it follows that $A\cup\{x,y\}\in\F$.
\end{itemize}
\end{lemma}

\subsubsection{Convex Geometries}
The way in which antimatroids encode the idea of ``convexity'' is expressed by the equivalent concept of an ``abstract convex geometry'', introduced by Edelman and Jamison \cite{edelman-jamison}.

Let $(E,\C)$ be a set system in which $\C$ is closed under intersections and $\emptyset\in\C$. This gives rise to a {\sf closure operator} $\tau:2^E\to 2^E$,
\begin{equation*}
\tau(X):=\bigcap\left\{ A\in\C, A\supseteq X\right\},
\end{equation*}
which satisfies the following properties:
\begin{itemize}
\item[] $\tau(\emptyset)=\emptyset$,
\item[] $A\subseteq\tau(A)$ for all $A\in 2^E$,
\item[] $A\subseteq B$ implies $\tau(A)\subseteq\tau(B)$,
\item[] $\tau(\tau(A))=\tau(A)$.
\end{itemize}
Conversely, every closure operator arises in this way. That is, if we are given a map $\tau$ satisfying the above properties, then $\C$ is the collection of sets satisfying $\tau(A)=A$---called the {\sf $\tau$-closed} sets. The triple $(E,\C,\tau)$ (where either $\C$ or $\tau$ is redundant) is called a closure space.

\begin{figure}
\begin{center}\input{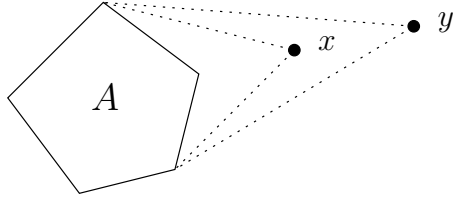}\end{center}
\caption{The anti-exchange axiom}
\label{fig:antiexchange}
\end{figure}

Closure operators are ubiquitous in combinatorics. If a closure space $(E,\C,\tau)$ satisfies the {\sf exchange axiom},
\begin{itemize}
\item[]
If $x,y\not\in\tau(A)$ then $x\in\tau(A\cup\{y\})$ implies $y\in\tau(A\cup\{x\})$,
\end{itemize}
then it is called a {\sf matroid}. In this case $\tau$ models the notion of ``linear span.'' If instead we wish to model the notion of ``convex hull'', we will require the following {\sf anti-exchange axiom},
\begin{itemize}
\item[]
If $x,y\not\in\tau(A)$ then $x\in\tau(A\cup\{y\})$ implies $y\not\in\tau(A\cup\{x\})$,
\end{itemize}
which is illustrated in Figure \ref{fig:antiexchange}.
In this case we say that $\tau$ is a {\sf convex closure} and that $(E,\C,\tau)$ is an {\sf abstract convex geometry}---or just a {\sf convex geometry}. Convex geometries and antimatroids are complementary structures.
\begin{lemma}{\bf \cite[Theorem III.1.3]{KLS}}
Given an accessible set system $(E,\F)$ on finite ground set $E$, let $\F^c=\{ E\setminus A:A\in\F\}$ denote the collection of complements of feasible sets. Then $(E,\F)$ is an antimatroid if and only if $(E,\F^c)$ is a convex geometry. In this case $\F$ consists of the open sets and $\F^c$ consists of the closed sets of a convex closure.
\end{lemma}

The contrast between the exchange and anti-exchange properties is one reason for the term {\em anti}-matroid.

The motivating example of a convex geometry is a pair $(E,\C)$ where $E\subseteq\R^n$ is a finite subset of Euclidean space and $\C$ is the collection of intersections of $E$ with convex subsets of $\R^n$.

\subsubsection{Join-Distributive Lattices}
\label{sec:JD}
Closure spaces, in turn, give rise to lattices. A partially-ordered set ({\sf poset}) $(P,\leq)$ is a finite set together with a reflexive, antisymmetric and transitive relation. If each pair of elements $x,y\in P$ possesses a least upper bound $x\vee y$ (their {\sf join}) and a greatest lower bound $x\wedge y$ (their {\sf meet}), then we call $(P,\leq)$ a {\sf lattice}. 

Given a finite closure space $(E,\C,\tau)$ let $(P_\C,\leq)$ denote the collection of $\tau$-closed sets, partially ordered by inclusion. Since $\C$ is closed under intersections, every pair $X,Y\in P_\C$ has a meet $X\wedge Y= X\cap Y$. Furthermore, since the ground set $E$ is $\tau$-closed, the collection of upper bounds $U(X,Y)=\{Z\in P_\C:X\subseteq Z \text{ and } Y\subseteq Z\}$ is nonempty and
\begin{equation*}
X\vee Y=\bigwedge_{Z\in U(X,Y)} Z
\end{equation*}
defines a join operation. Hence $(P_\C,\leq)$ is a lattice.

The relationship between closure spaces and lattices has been fruitful and the two most important types of closure spaces---matroids and antimatroids---have been classified in terms of their lattice structure. To express this we need some poset notation.

Given elements $x,y$ in a poset $(P,\leq)$, we say that $y$ {\sf covers} $x$ (and write $x\prec y$) when $x\leq y$ and there does not exist $z\in P$ such that $x< z< y$. If $P$ has a minimum element $\hat{0}$ then the elements covering $\hat{0}$ are called {\sf atoms}. More generally, the atoms of an interval $[x,y]=\{z\in P:x\leq z\leq y\}$ are the elements in $[x,y]$ that cover $x$. We call a lattice $P$ {\sf upper-semimodular} if whenever $x$ and $y$ both cover $x\wedge y$ it follows that $x\vee y$ covers both $x$ and $y$. A lattice $P$ is called {\sf geometric} if it is upper-semimodular and if every element in $P$ can be written as a join of atoms. The following characterization was known to Birkhoff, Whitney and MacLane in the 1930s (see \cite[Page 179]{gratzer}).

\begin{theorem}{\bf \cite[Theorem I.1.7]{KLS}}
A finite lattice is geometric if and only if it occurs as the lattice of closed sets of an exchange closure.
\end{theorem}

The characterization of anti-exchange closures is more recent and is due to Edelman. We say that a lattice is {\sf boolean} if it is isomorphic to the collection of all subsets of a finite set under inclusion. Given an element $x$ in a lattice and a set $Y$ of elements covering $x$, we call the interval
$$ \left[ x, \bigvee_{y\in Y} y\right] $$
an {\sf atomic interval}. If $\hat{1}$ is the maximum element in a lattice, then an element $x\neq\hat{1}$ is  called {\sf meet-irreducible} if $x=y\wedge z$ implies either $x=y$ or $x=z$. A lattice that satisfies any of the following equivalent conditions is called {\sf join-distributive}.
\begin{theorem}{\bf \cite[Theorems 1.7 and 1.9]{AGT}}
Given a finite lattice $(P,\leq)$, the following statements are equivalent.
\begin{enumerate}
\item Every atomic interval in $P$ is boolean.
\item Every element of $P$ has a {\em unique} irredundant decomposition as a meet of meet-irreducible elements.
\label{item:dilworth}
\item $P$ is upper-semimodular and it satisties the {\sf meet-semidistributive property}: for all $x,y,z\in P$, we have
\begin{equation*}
x\wedge y=x\wedge z\quad\Longrightarrow\quad x\wedge y=x\wedge(y\vee z).
\end{equation*}
\end{enumerate}
\end{theorem}

These lattices were first considered by Dilworth \cite{dilworth}, for whom condition \ref{item:dilworth} was the defining property. Edelman's theorem is the following.

\begin{theorem}{\bf \cite[Theorem 3.3]{edelman:meet}}
\label{th:edelman}
A finite lattice is join-distributive if and only if it occurs as the lattice of {\em open} sets of an anti-exchange closure.\footnote{Edelman \cite{edelman:meet} used the term {\sf meet-distributive} for the lattice of {\em closed} sets.}
\end{theorem}

That is, a finite join-distributive lattice is precisely the lattice of feasible sets of some antimatroid. In this sense, geometric lattices and join-distributive lattices are opposite. In fact, since a geometric lattice is itself an atomic interval, the only lattices that are both geometric and join-distributive are the boolean lattices.

\subsection{Supersolvable Join-Distributive Lattices}

In \cite{stanley:ss} Stanley introduced another generalization of distributive lattices---the {\sf supersolvable lattices}. (See also \cite[Example 3.13.4]{stanley:ec1}.) The motivating example of a supersolvable lattice is the lattice of subgroups of a supersolvable group. In the same paper Stanley also discussed join-distributive lattices---under the name $1$-lattices---and he used the code {\sf $1$-SS} to refer to lattices that are both supersolvable and join distributive. Hawrylycz and Reiner \cite{HR} later found an important class of examples: the lattice of closure relations on any finite poset is {\sf $1$-SS}.

In this paper we will explore a new source of examples. Hence we find it convenient to give a characterization of {\sf $1$-SS} lattices in the spirit of Edelman's  theorem (Theorem \ref{th:edelman}). First let us review the notion of ``supersolvability.''

A finite lattice $(P,\leq)$ is called {\sf distributive} if either of the equivalent ``distributive laws'',
\begin{itemize}
\item[] $x\wedge(y\vee z)=(x\wedge y)\vee(x\wedge z)$,
\item[] $x\vee(y\wedge z)=(x\vee y)\wedge(x\vee z)$,
\end{itemize}
holds for all $x,y,z\in P$.
Given a subset $X\subseteq P$, the smallest sublattice of $P$ containing $X$ is called the sublattice {\sf generated by} $X$. A maximal chain $\mathfrak{m}$ in $P$ is called an {\sf M-chain} if together with any other maximal chain in $P$ it generates a distributive sublattice. If such an {\sf M}-chain exists we say that $P$ is a {\sf supersolvable lattice}.

A characterizing property of supersolvable latttices is the fact that they have a certain kind of ``edge-labelling.'' Given a poset $P$, its {\sf Hasse diagram} is a directed graph with an edge from $x$ to $y$ for each cover relation $x\prec y$ in $P$. (The edges are directed ``up'' in the diagram.) An {\sf edge-labelling} of $P$ is just a labelling of the edges of its Hasse diagram by the elements of some ordered set---typically the integers. Now suppose that $P$ is supersolvable with {\sf M}-chain
\begin{equation*}
\hat{0}=a_0\prec a_1\prec\cdots \prec a_n=\hat{1},
\end{equation*}
and define an edge-labelling $\lambda$ by setting
\begin{equation}
\label{eq:edgelabel}
\lambda(x,y):=\min\left\{ i: y\leq a_i\vee x\right\}
\end{equation}
for each cover $x\prec y$ in $P$. Stanley \cite{stanley:ss} showed that this labelling has the following properties:
\begin{itemize}
\item[] The labels on any maximal chain in $P$ form a permutation of the index set, which we may assume is $\{1,\ldots,n\}$. It follows that the labels on the maximal chains in an interval $[x,y]\subseteq P$ are permutations of some fixed subset $I_{[x,y]}\subseteq \{1,\ldots,n\}$.
\item[] In any interval $[x,y]\subseteq P$ there exists a unique maximal chain labelled by the increasing permutation of $I_{[x,y]}$.
\end{itemize}
Such a labelling is called an $S_n$ {\sf EL-labelling} of $P$ (where {\rm EL} stands for ``edge-lexicographic''). McNamara proved that the existence of an $S_n$ \EL-labelling characterizes a supersolvable lattice.

\begin{lemma}{\bf \cite[Theorem 1]{mcnamara}}
\label{lem:mcnamara}
A finite lattice is supersolvable if and only if it possesses an $S_n$ \EL-labelling. In particular, the unique increasing maximal chain in this labelling is an {\sf M}-chain.
\end{lemma}

Now we can give our characterization of supersolvable join-distributive lattices, based on the following definition. We say that $\leq_E$ is a {\sf total order} on $E$ if either $x\leq_E y$ or $y\leq_E x$ for all $x,y\in E$.

\begin{definition}
\label{def:ssantimatroid}
A set system $(E,\F,\leq_E)$, with a total order $\leq_E$ on the ground set $E$, is called a {\sf supersolvable antimatroid} if it satisfies $\emptyset\in\F$ and the following condition:
\begin{itemize}
\item[\eq\label{eq:ssantimatroid}]
Given feasible sets $A,B\in\F$ with $B\not\subseteq A$ and $x=\min_{\leq_E}(B\setminus A)$, it follows that $A\cup\{x\}\in\F$.
\end{itemize}
\end{definition}

In particular, a supersolvable antimatroid is an antimatroid. Indeed, taking $A=\emptyset$, condition \eqref{eq:ssantimatroid} implies that every prefix of a feasible set---with respect to the order on $E$---is feasible. Thus the maximum element of any feasible set may deleted, and we conclude that $(E,\F)$ is an accessible system. Finally, note that condition \eqref{eq:ssantimatroid} strengthens condition \eqref{eq:antimatroid}.

\begin{theorem}
\label{th:ssjd}
Let $(E,\F)$ be an antimatroid with join-distributive lattice $P$ of feasible sets. The following statements are equivalent.
\begin{enumerate}
\item There exists a total order on $E$ with respect to which $(E,\F)$ is a supersolvable antimatroid.
\label{ssam1}
\item $P$ is a supersolvable lattice.
\label{ssam2}
\item There exists an order on $E$ with respect to which the natural edge-labelling of $P$ by $E$ is an $S_n$ \EL-labelling.
\label{ssam3}
\end{enumerate}
\end{theorem}

\begin{proof}
We will show that $\ref{ssam1}\Rightarrow\ref{ssam2}\Rightarrow\ref{ssam3}\Rightarrow\ref{ssam1}$.

First, suppose that $(E,\F)$ is supersolvable with respect to a certain total order $(E,\leq_E)$ and consider the join-distributive lattice $(P,\leq)$ of feasible sets. Note that a cover relation $X\prec Y$ in $P$ is naturally labelled by the element $x\in E$ where $\{x\}=Y\setminus X$. We claim that this is an $S_n$ \EL-labelling with respect to the order $\leq_E$. Indeed, consider an interval $[A,B]$ in $P$. By property \eqref{eq:antimatroid} of antimatroids we know that every maximal chain in $[A,B]$ is labelled by some permutation of the set $B\setminus A\subseteq E$, and condition \eqref{eq:ssantimatroid} guarantees that the unique increasing permutation occurs. Thus $P$ possessees an $S_n$ \EL-labelling and by Lemma \ref{lem:mcnamara} it is supersolvable.

Next, suppose that $P$ is supersolvable with {\sf M}-chain
\begin{equation*}
\emptyset=A_0\prec A_1\prec\cdots\prec A_n=E,
\end{equation*}
and set $\{x_i\}=A_i\setminus A_{i-1}$ for all $1\leq i\leq n$. This defines a total order $x_1<_E x_2<_E\cdots<_E x_n$ on the set $E$. Now consider a cover relation $X\prec Y$ in $P$. Since the join operation in $P$ is just union of feasible sets, we conclude that the edge-labelling \eqref{eq:edgelabel} satisfies $\lambda(X,Y)=i$, where $Y\setminus X=\{x_i\}$. That is, the labelling $\lambda$ coincides with the natural edge-labelling by $E$. Since $\lambda$ is an $S_n$ EL-labelling, so is the edge-labelling by $E$.

Finally, let $(E,\leq_E)$ be a total order such that the edge-labelling of $P$ by $E$ is $S_n$ \EL. Now consider $A$ and $B$ in $\F$ (hence also $A\cup B\in\F$) with $B\not\subseteq A$. By assumption there is a maximal chain in the interval $[A,A\cup B]$ that is labelled by the unique increasing permutation of $(A\cup B)\setminus A=B\setminus A\subseteq E$. If $x\in E$ is the first label on the chain, we conclude that $x=\min_{\leq_E}(B\setminus A)$ and $A\cup\{x\}\in\F$. Hence $(E,\F)$ is a supersolvable antimatroid.
\end{proof}

\begin{figure}
\begin{center}\input{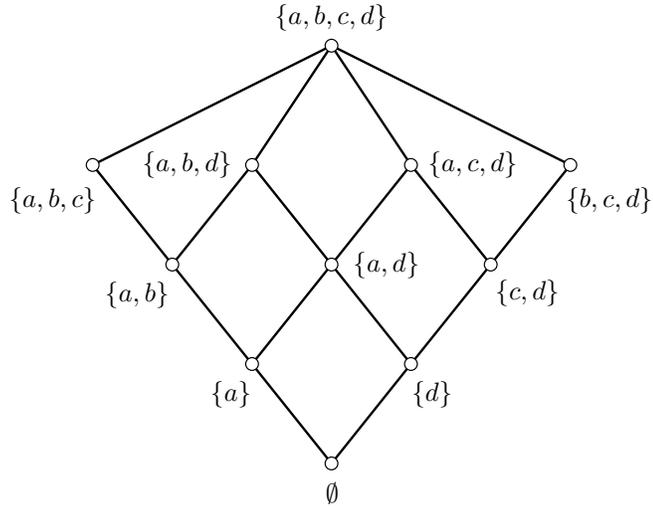}\end{center}
\caption{A join-distributive lattice that is not supersolvable}
\label{fig:jdnotss}
\end{figure}

Both ``supersolvable'' and ``join-distributive'' are generalizations of the concept ``distributive.'' However, they are distinct concepts. For instance, the lattice
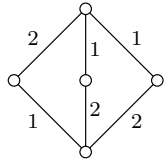
\begin{center}\input{M3.pstex_t}\end{center}
is supersolvable (it has an $S_n$ EL-labelling) but it is not join-distributive since the lattice itself is an atomic interval that is not boolean. The next example illustrates that a join-distributive lattice need not be supersolvable.

\begin{example}
Consider the set $E=\{a<b<c<d\}$ of four distinct points in the line $\R$. Recall that the intersections of $E$ with complements of convex sets in $\R$ are the feasible sets of an antimatroid, and the inclusion order on these is a join-distributive lattice. We display this lattice in Figure \ref{fig:jdnotss}. It is easy to see, however, that this lattice is not supersolvable. For instance, the intervals $[\emptyset,\{a,b,c\}]$ and $[\emptyset,\{b,c,d\}]$ are chains, with edge-label sequences $(a,b,c)$ and $(d,c,b)$. It is impossible to order the set $E$ in such a way that both of these sequences are increasing.
\end{example}

At the moment we know relatively little about supersolvable antimatroids in general. This may be an interesting avenue for further study. In the next sections we will describe a natural class of supersolvable antimatroids arising from Coxeter groups.

\section{Sorting in a Coxeter Group}
\label{sec:sorting}

In the next two sections let $(W,S)$ be an arbitrary Coxeter system and consider a word $\om\in S^*$ in the generators $S$. At first we suppose that $\om$ is an arbitrary finite word. However, we will see below (Corollary \ref{cor:reduced}) that we lose nothing by taking $\om$ to be reduced. In Section \ref{sec:infinite} we will allow $\om$ to have infinite length.

We call $\om$ the {\sf sorting word} and we will use it to ``sort'' each element of $W$. This leads to a notion of ``$\om$-sorted'' words and an associated partial order on the group.

\subsection{Sorted Words}
First we set down some notation. Let $\om=(\om_1,\om_2,\ldots,\om_m)\in S^*$ denote the sorting word. We will typically identify $\om$ with the {\sf ground set} $I(\om):=\{1,2,\ldots,m\}$. We say that $\alpha=(\alpha_1,\alpha_2,\ldots,\alpha_k)\in S^*$ is a {\sf subword} of $\om$ (and we write $\alpha\subseteq\om$) if we have
\begin{equation*}
(\alpha_1,\alpha_2,\ldots,\alpha_k)=(\om_{i_1},\om_{i_2},\ldots,\om_{i_k})
\end{equation*}
for some $1\leq i_1<i_2<\cdots <i_k\leq m$. In this case, the {\sf index set} of $\alpha$ is $I(\alpha):=\{i_1,i_2,\ldots,i_k\}\subseteq I(\om)$. Thus $I(\alpha)\subseteq I(\beta)\subseteq I(\om)$ implies that $\alpha$ is a subword of $\beta$; the converse is not true.

\begin{warning}
Because the sorting word $\om$ may contain repeated letters, the index set $I(\alpha)$ of a subword $\alpha\subseteq\om$ may not be uniquely recoverable from $\alpha$. Thus we will always identify a subword with the pair $(\alpha,I(\alpha))$. (We may omit mention of the index set $I(\alpha)$ when no confusion will result.) Two subwords $\alpha,\beta\subseteq\om$ are {\sf equal} when $I(\alpha)=I(\beta)$. We will write $\alpha\cup\beta$ and $\alpha\cap\beta$ to denote the words corresponding to index sets $I(\alpha)\cup I(\beta)$ and $I(\alpha)\cap I(\beta)$, respectively.
\end{warning}

To each word $\alpha=(\alpha_1,\alpha_2,\ldots,\alpha_k)\in S^*$ we associate the group element
\begin{equation*}
\langle\alpha\rangle:=\alpha_1\alpha_2\cdots\alpha_k\in W.
\end{equation*}
Note that the correspondence $\alpha\mapsto\langle\alpha\rangle$ is not injective. We say that $\alpha$ is a {\sf reduced word} (for $\langle\alpha\rangle$) if there does not exist another word for $\langle\alpha\rangle$ of shorter length. In this case, the group element $\langle\alpha\rangle\in W$ has {\sf length} $k$, and we write $\ell(\langle\alpha\rangle)=k$. Notice that $\ell(w)=\ell(w^{-1})$ for all $w\in W$ since we may reverse a reduced word for $w$ to obtain one for $w^{-1}$.

In this paper we will consider the collection $2^{I(\om)}$ of subsets of the ground set $I(\om)=\{1,2,\ldots,m\}$ together with the total {\sf lexicographic order} $\leq_{\lex}$. Given $A$ and $B$ subsets of $I(\om)$, we will say that $A\leq_\lex B$ if either $A=B$ or the minimum element of $(A\cup B)\setminus (A\cap B)$ is in $A$.

We can now define sorted words.

\begin{definition}
A subword $\alpha\subseteq\om$ of the sorting word is called {\sf $\om$-sorted} if
\begin{enumerate}
\item $\alpha$ is a reduced word,
\item $I(\alpha)=\min_{\leq_\lex}\left\{ I(\beta)\subseteq I(\om): \text{$\beta$ is reduced and }\langle\beta\rangle=\langle\alpha\rangle\right\}$.
\end{enumerate} 
\end{definition}
That is, $\alpha$ is $\om$-sorted if it is the lexicographically-least reduced word for $\langle\alpha\rangle$ among subwords of $\om$.

\subsection{The Sorting Algorithm}

Let $W_\om\subseteq W$ denote the subset of group elements that occur as subwords of the sorting word $\om$. Each element of $W_\om$ corresponds to a unique $\om$-sorted word, and we may think of this as a canonical form for the group element.

Recognizing canonical forms is a fundamental problem. In this section we given an algorithmic characterization of $\om$-sorted words.

\begin{sortingalg}
Given a group element $u \in W_\om$, let $\alpha\subseteq\om=(\om_1,\ldots,\om_m)$ be a subword of the sorting word such that $u=\langle\alpha\rangle$. We define another subword $\omsort(\alpha)\subseteq\om$ as follows. Begin by setting $X:=\emptyset$ and $x:=u$. For $i$ from $1$ to $m$ do
\begin{itemize}
\item If $\ell(\om_i x)=\ell(x)-1$ then put $X:=X\cup\{i\}$ and $x:=\om_i x$.
\item If $\ell(\om_i x)=\ell(x)+1$ do nothing.
\end{itemize}
Let $\omsort(\alpha)=\omsort(u)$ be the subword of $\om$ with index set $X$. This is called the {\sf $\om$-sorted word} for $\alpha\in S^*$ and $u\in W$.
\end{sortingalg}

That is, we proceed through the entries of $\om$, checking successively whether $\om_i\in S$ is a {\sf left descent} of our group element $x$. If it is, then we record the index $i$ and replace our group element $x$ by $\om_i x$. Otherwise, we do nothing. We illustrate this algorithm with an example.

\begin{example}
\label{ex:sort}
Consider a Coxeter system $(W,S)$ of type $A_{n-1}$. That is, let $W=\mathfrak{S}_n$ be the group of permutations of $\{1,2,\ldots,n\}$ with generating set $S$ of adjacent transpositions,
\begin{equation*}
S=\left\{ s_i=(i,i+1) : 1\leq i\leq n-1\right\}.
\end{equation*}
We will express a permutation $\sigma\in\mathfrak{S}_n$ using the one-line notation 
\begin{equation*}
\sigma(1)\sigma(2)\cdots\sigma(n).
\end{equation*}
Notice that $s_i$ is a left descent for $\sigma\in\mathfrak{S}_n$ (that is, $\ell(s_i\sigma)=\ell(\sigma)-1$) precisely when the symbols $i$ and $i+1$ are out of order in the word for $\sigma$ (that is, when $\sigma^{-1}(i+1)<\sigma^{-1}(i)$).

Now let $n=5$ and fix the sorting word
\begin{equation*}
\om=(\om_1,\ldots,\om_{10})=(s_1,s_2,s_3,s_4,s_3,s_2,s_1,s_2,s_3,s_2).
\end{equation*}
(Incidentally, this is a reduced word for the longest element $w_\circ=54321\in\mathfrak{S}_5$ so that $W_\om$ is the full group. This follows from the fact that $w_\circ$ is the unique maximum element under Bruhat order --- see Section \ref{sec:weakBruhat}.) We compute the $\om$-sorted word for $\sigma=41532$ in Table \ref{table:omsort}.
\begin{table}
\begin{center}
\begin{tabular}{c|c|c|c|c}
step & reflection & $x\in W$ & descent? & index set $X$\\
\hline
1 & $s_1=(12)$ & $41532$ & no  & $\{$\\ 
2 & $s_2=(23)$ & $41532$ & yes & $\{2$\\ 
3 & $s_3=(34)$ & $41523$ & yes & $\{2,3$\\
4 & $s_4=(45)$ & $31524$ & yes & $\{2,3,4$\\
5 & $s_3=(34)$ & $31425$ & no   & $\{2,3,4$\\
6 & $s_2=(23)$ & $31425$ & yes & $\{2,3,4,6$\\
7 & $s_1=(12)$ & $21435$ & yes & $\{2,3,4,6,7$\\
8 & $s_2=(23)$ & $12435$ & no   & $\{2,3,4,6,7$\\
9 & $s_3=(34)$ & $12435$ & yes & $\{2,3,4,6,7,9$\\
10&$s_2=(23)$ & $12345$ & no  & $\{2,3,4,6,7,9\}$
\end{tabular}
\end{center}
\caption{An example of the $\om$-sorting algorithm}
\label{table:omsort}
\end{table}
Since the resulting index set is $\{2,3,4,6,7,9\}$, we obtain the $\om$-sorted word
\begin{equation*}
\omsort(41532)=\omsort(\sigma)=(s_2,s_3,s_4,s_2,s_1,s_3).
\end{equation*}
\end{example}

One should imagine that the algorithm converts a group element $u\in W_\om$ to the identity element and that the word $\omsort(u)$ records the steps in this process.

The sorting algorithm has the following properties.

\begin{lemma}
\label{lemma:algorithm}
Let $\alpha\subseteq\om$ be a subword of the sorting word. We have
\begin{enumerate}
\item $\omsort(\alpha)$ is a reduced word for $\langle\alpha\rangle\in W$.
\label{item:sort1}
\item $\alpha$ is $\om$-sorted if and only if $\omsort(\alpha)=\alpha$.
\label{item:sort2}
\end{enumerate}
\end{lemma}

\begin{proof}
\ref{item:sort1}. Suppose that $\omsort(\alpha)=(\gamma_1,\ldots,\gamma_k)$, so the algorithm will terminate with $x=\gamma_k\cdots\gamma_2\gamma_1\langle\alpha\rangle$. We also have
\begin{equation*}
\ell(x)=\ell(\langle\alpha\rangle)-k,
\end{equation*}
since by construction each multiplication with a generator decreases the length of $x$ by $1$. Thus we will be done if we can show that the algorithm terminates with $x=1\in W$---the identity element.

To show this we use induction and the Exchange Property (see Section \ref{sec:weakBruhat} below)---hereafter invoked just as ``Exchange.'' Suppose that the following statement holds for some index $i$:
\begin{quote}
We have completed $i-1$ steps of the algorithm and we currently have $x=\langle\alpha'\rangle$, where $\alpha'=(\alpha'_1,\ldots,\alpha'_\ell)$ is a subword of $(\om_i,\om_{i+1},\ldots,\om_m)$.
\end{quote}
Note that this statement is true with $i=1$ and $\alpha'=\alpha$. On the $i$th step, we perform the multiplication $\om_i\langle\alpha'\rangle$. If $\ell(\om_i\langle\alpha'\rangle)=\ell(\langle\alpha'\rangle)+1$, then $\alpha'$ is in fact a subword of $(\om_{i+1},\ldots,\om_m)$. Thus we do nothing, and the statement remains true. On the other hand, suppose that $\ell(\om_i\langle\alpha'\rangle)=\ell(\langle\alpha'\rangle)-1$. In this case, by Exchange there exists $1\leq j\leq \ell$ such that
\begin{equation*}
\om_i\langle\alpha'\rangle=\om_i\alpha'_1\cdots\alpha'_\ell=\alpha'_1\cdots\hat{\alpha}'_j\cdots\alpha'_\ell.
\end{equation*}
Replacing $\alpha'$ by $(\alpha'_1,\ldots,\hat{\alpha}'_j,\ldots,\alpha'_\ell)$, the statment remains true. By induction, the statment is true for $i=m$---thus the algorithm terminates with $\alpha'$ equal to the empty word as desired.

\ref{item:sort2}. First suppose that $\alpha\subseteq\om$ is {\em not} $\om$-sorted and let $\beta=\omsort(\alpha)$. If $\alpha=(\alpha_1,\ldots,\alpha_k)$ and $\beta=(\beta_1,\ldots,\beta_\ell)$,  let $j$ be the minimum integer such that $\alpha_j\neq\beta_j$. Since $I(\beta)<_\lex I(\alpha)$, this implies that $\beta_j=\om_{j'}$ with $\om_{j'}\not\in\alpha$. Now let us apply the $\om$-sorting algorithm to $\alpha$. At the $(j'-1)$th step we will have
\begin{equation*}
x=\alpha_j\alpha_{j+1}\cdots\alpha_k=\beta_j\beta_{j+1}\cdots\beta_\ell.
\end{equation*}
Next, applying $\om_{j'}$ on the left we get
\begin{equation*}
\ell(\om_{j'}x)=\ell(\beta_{j+1}\cdots\beta_\ell)=\ell(x)-1,
\end{equation*}
which implies that $\om_{j'}\in\omsort(\alpha)$. Since $\om_{j'}\not\in\alpha$, we conclude that $\omsort(\alpha)$ is not equal to $\alpha$.

Conversely, suppose that $\omsort(\alpha)\neq\alpha$. Since $\alpha$ is a word for $\langle\alpha\rangle$ and $\omsort(\alpha)$ is a {\em reduced} word for $\langle\alpha\rangle$ (by part \ref{item:sort1}), there exists a minimum integer $j$ such that $\om_j\in\omsort(\alpha)$ and $\om_j\not\in\alpha$. In this case,
\begin{equation*}
\omsort(\alpha)=(\alpha\cap (\om_1,\ldots,\om_{j-1}))\cup(\omsort(\alpha)\cap (\om_j,\ldots,\om_m))
\end{equation*}
is a reduced word for $\langle\alpha\rangle$ that is strictly lexicographically-less than $\alpha$. Hence $\alpha$ is not $\om$-sorted.
\end{proof}

This result gives us a convenient way to recognize $\om$-sorted words, which we will use in later proofs.

\begin{corollary}
\label{cor:sorted}
A word $\alpha\subseteq\om=(\om_1,\ldots,\om_m)$ is $\om$-sorted if and only if there does {\em not} exist $1\leq j< m$ with $\om_j\not\in\alpha$ such that $\om_j$ is a left descent of the group element $\langle\alpha\cap(\om_{j+1},\ldots,\om_m)\rangle\in W$.
\end{corollary}

\subsection{Remark---Classical Sorting Algorithms}
We have two justifications for our use of the term ``sorting.''

First, Reading \cite{reading:csort1} has defined the notion of ``Coxeter-sorting'' for elements in a Coxeter group (see Section \ref{sec:cyclic} below). He notes that his Coxeter-sorting algorithm is related to the classical ``stack-sorting'' algorithm, described by Knuth \cite[Exercise 2.2.1.4--5]{knuth:taocp}. Reading's work is the main motivation behind the current paper.

More generally, Knuth \cite[Chapter 8]{knuth:aah} has described a framework for a wide variety of classical sorting algorithms. He defines a {\sf comparator} $[i:j]$ as a map that operates on a sequence of numbers $(x_1,\ldots,x_n)$, replacing $x_i$ and $x_j$ respectively by $\min(x_i,x_j)$ and $\max(x_i,x_j)$. A {\sf sorting network} is a sequence of comparators that will sort any given sequence $(x_1,\ldots,x_n)$, and a {\sf primitive sorting network} consists entirely of comparators of the form $[i:i+1]$. Thus, our sorting word $\om$ and $\om$-sorting algorithm may be thought of as a generalized ``primitive sorting network'' for the set of ``sequences'' $W_\om\subseteq W$.

Furthermore, Knuth notes that one may restrict attention to the {\em irredundant} primitive sorting networks, which correspond to commutation classes of reduced words for the longest permutation $w_\circ=n(n-1)\cdots 321$. (In Section \ref{sec:ssantimatroid} we will show that the same reduction can be made in general.) Various reduced words for $w_\circ$ then correspond to different classical sorting algorithms.
\begin{figure}
\begin{center}\input{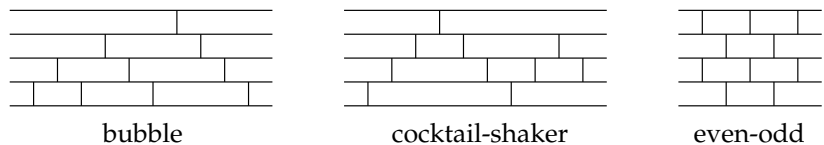}\end{center}
\caption{Some classical sorting algorithms}
\label{fig:knuth}
\end{figure}
Figure \ref{fig:knuth} (adapted from Knuth \cite[page 29]{knuth:aah}, which contains an error) shows some reduced words for $w_\circ=54321$ in $\mathfrak{S}_5$ corresponding to three classical sorting algorithms. Here, ``bubblesort'' corresponds to the lexicographically first reduced word. In our Example \ref{ex:sort}, we have performed the ``cocktail-shaker sort.''

\subsection{The Sorting Order}
To end this section we define a natural partial order on the collection of $\om$-sorted words, or equivalently on the set of group elements $W_\om\subseteq W$.

\begin{definition}
Given a Coxeter system $(W,S)$ and a sorting word $\om\in S^*$ let $(P_\om,\leq_\om)$ denote the set of $\om$-sorted words together with the subword inclusion order,
\begin{equation*}
\alpha\leq_\om\beta \quad\Longleftrightarrow\quad I(\alpha)\subseteq I(\beta).
\end{equation*}
We also define a partial order on $W_\om$ by identifying a group element with its $\om$-sorted word,
\begin{equation*}
u \leq_\om w \quad\Longleftrightarrow\quad \omsort(u)\,\leq_\om\,\omsort(w).
\end{equation*}
\end{definition}

In the next section we will see that the sorting order has many remarkable properties and it is closely related to other important orders on the group $W$.

\section{Properties of the Sorting Order}
\label{sec:properties}

\subsection{Between Weak and Bruhat}
\label{sec:weakBruhat}
Let $(W,S)$ be an arbitrary Coxeter system. Recall that $(s_1,\ldots,s_k)\in S^*$ is called a {\sf reduced word} for the group element $w\in W$ when $w=s_1\cdots s_k$ and there does not exist a word for $w$ of length less than $k$. In this case $\ell(w)=k$ is the {\sf length} of the element $w\in W$. There are two classical and important partial orders on the group $W$, both based on the combinatorics of reduced words (this topic is covered thoroughly in \cite{bjorner-brenti}).

\begin{weakorder}
Given $u,w\in W$, we write $u\leq_\RW w$ if $u$ occurs as a ``prefix'' of $w$---that is, if there exists a reduced word $w=s_1s_2\cdots s_\ell$ for $w$ and an integer $1\leq k\leq \ell$ such that $u=s_1s_2\cdots s_k$ is a reduced word for $u$. This is called the {\sf weak order} on $W$.
\end{weakorder}

We have actually defined the ``right'' weak order. There is a corresponding ``left'' weak order and the two are isomorphic via the map $w\mapsto w^{-1}$ which exchanges prefixes and suffixes. The Hasse diagram of the weak order is just the right Cayley graph of $W$ with respect to the generating set $S$. In general the weak order is graded by the length function $\ell:W\to\Z$ and it is a meet-semilattice---it possesses meets, but not joins. However, when $W$ is finite there exists a unique element $w_\circ\in W$ of maximum length---called the ``longest element''---which satisfies $w\leq_\RW w_\circ$ for all $w\in W$. Hence the weak order on a finite Coxeter group is a lattice.

\begin{bruhatorder}
Given $u,w\in W$, we write $u\leq_\B w$ if $u$ occurs as a ``subword'' of $w$---that is, if there exists a reduced word $w=s_1s_2\cdots s_\ell$ for $w$ and integers $1\leq i_1<\cdots <i_k\leq\ell$ such that $u=s_{i_1}\cdots s_{i_k}$ is a reduced word for $u$. This is called the {\sf Bruhat order} on $W$.
\end{bruhatorder}

(Moreover, if $u\leq_B w$, it turns out that $u$ occurs as a subword of {\bf any} reduced word for $w$.) It is easy to see from the definition that the Bruhat order is also a graded poset, ranked by the length function. Since a ``prefix'' is a ``subword'' we note that $u\leq_\RW w$ implies $u\leq_\B w$ for all $u,w\in W$---that is, Bruhat order is a {\sf poset extension} of the weak order. However Bruhat order is neither a meet- nor a join-semilattice. The Bruhat order arises in applications as the inclusion order on closures of Schubert cells in the generalized flag variety corresponding to $W$.

Since the Bruhat order is an extension of the weak order with the same rank function, one may obtain the Hasse diagram of Bruhat order from the hasse diagram of weak order by adding some extra edges (cover relations). Let $T=\{ wsw^{-1}: w\in W, s\in S\}$ denote the generating set of {\sf reflections}. By definition, the cover relations in weak order have the form $u\prec w$ where $w=us$ for some $s\in S$. It is also true---but not obvious---that the covers in Bruhat order have the form $u\prec w$ where $w=ut$ for some reflection $t\in T$ such that $\ell(w)=\ell(u)+1$. We will find that the sorting orders are intermediate between weak and Bruhat order since they include some but not all of these extra covers.

A good first example is the dihedral group.

\begin{figure}
\begin{center}\input{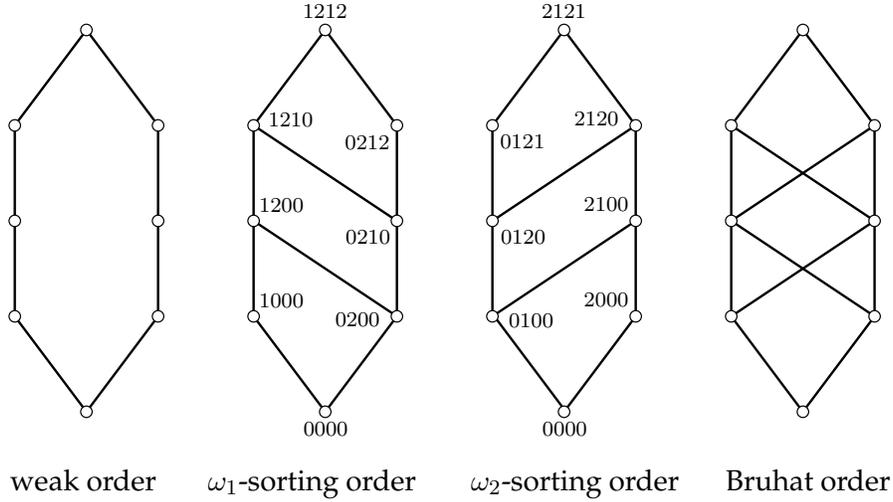}\end{center}
\caption{Comparison of weak, Bruhat and sorting orders}
\label{fig:dihedrals}
\end{figure}

\begin{example}
\label{ex:dihedrals}
Let $W=I_2(m)$ be the dihedral group of order $2m$ with Coxeter generators $S=\{s_1,s_2\}$. In this case the longest element $w_\circ\in W$ is of length $m$, and it has exactly two reduced words: $\om_1:=(s_1,s_2,s_1,s_2,\ldots)$ and $\om_2:=(s_2,s_1,s_2,s_1,\ldots)$. The following is a convenient notation for $\om_1$- and $\om_2$-sorted words: For example, when $m=4$ and $\om_1=(s_1,s_2,s_1,s_2)$, the subword $\alpha=(s_2,s_1,s_2)$ is $\om_1$-sorted with index set $I(\alpha)=\{2,3,4\}\subseteq I(\om_1)=\{1,2,3,4\}$. We encode both the word and the index set simultaneously with the string $0212$, where the zeroes are placeholders. Figure \ref{fig:dihedrals} displays the Hasse diagrams of the weak order, Bruhat order and both sorting orders on the group $I_2(4)$ (a.k.a. $B_2$).
\end{example}

Notice that the sorting orders on the full group $I_2(m)$ occur ``between'' the weak and Bruhat orders. To prove that this property holds in general, we will need the following well-known lemma regarding the combinatorics of reduced words (see \cite[Theorem 1.4.3]{bjorner-brenti} for a proof).

\begin{exchange} Let $(s_1,s_2,\ldots,s_k)\in S^*$ be a reduced word for $w\in W$ and suppose that $\ell(tw)<\ell(w)$ for some $t\in T$. Then there exists $1\leq i\leq k$ such that
\begin{equation*}
tw=s_1\cdots \hat{s_i}\cdots s_k.
\end{equation*}
\end{exchange}

Here the notation $\hat{s_i}$ indicates that the symbol $s_i$ has been deleted from the word. Note that this new word for $tw$ is not necessarily reduced since we might have $\ell(tw)<\ell(w)-1$. A version of Exchange also holds when we multiply on the right by a reflection $t\in T$---we see this by applying Exchange to $w^{-1}$ and noting that $\ell(tw^{-1})=\ell((wt)^{-1})=\ell(wt)$.

\begin{theorem}
\label{th:weakBruhat}
Given a Coxeter system $(W,S)$ and a sorting word $\om\in S^*$, the $\om$-sorting order extends the weak order on $W_\om$ and is extended by the Bruhat order on $W_\om$. That is, for all $u,w\in W_\om$ we have
\begin{equation*}
u\leq_\RW w\quad\Longrightarrow\quad u\leq_\om w \quad\Longrightarrow\quad u\leq_\B w.
\end{equation*}
\end{theorem}

\begin{proof}
The fact that Bruhat order extends $\om$-sorting order is immediate since the $\om$-sorting order is a special case of subword inclusion.

Now suppose that $u\leq_\RW w$---that is, there exists a reduced word $\beta=(s_1,\ldots,s_\ell)$ for $w$ and an integer $1\leq k\leq\ell$ such that $\alpha=(s_1,\ldots,s_k)$ is a reduced word for $u$. To demonstrate that $u\leq_\om w$ we must show that $\omsort(\alpha)$ is a subword of $\omsort(\beta)$. We do this by induction.

Let the sorting word be $\om=(\om_1,\om_2,\ldots,\om_m)\in S^*$ and perform the $\om$-sorting algorithm on $\alpha$ and $\beta$. In the first step we multiplity on the left by $\om_1\in S$. There are three cases:
\begin{enumerate}
\item If $\om_1$ is a left descent of $u$ (that is, if $\ell(\om_1 u)=\ell(u)-1$), then, by the Exchange property, there exists an integer $1\leq i\leq k$ such that $\alpha':=(s_1,\ldots,\hat{s_i},\ldots,s_k)$ is a reduced word for $su$. In this case we see that $s$ is a left descent of $w$ and $\beta':=(s_1,\ldots,\hat{s_i},\ldots,s_\ell)$ is a reduced word for $sw$. We find that $\om_1$ is in both $\omsort(\alpha)$ and $\omsort(\beta)$.
\item If $\om_1$ is a left descent of $w$, but not of $u$, then by Exchange there exists an integer $1\leq i\leq \ell$ such that $\beta':=(s_1,\ldots,\hat{s_i},\ldots,s_\ell)$ is a reduced word for $sw$. Since $\om_1$ is not a left descent of $u$ we must have $k<i$, hence we set $\alpha':=\alpha$. In this case $\om_1$ is in $\omsort(\beta)$ but not in $\omsort(\alpha)$.
\item Finally, if $\om_1$ is not a left descent of $w$ then we set $\alpha':=\alpha$ and $\beta':=\beta$. This time $\om_1$ is in neither $\omsort(\alpha)$ nor $\omsort(\beta)$.
\end{enumerate}
In any case we find that $\om_1\in\omsort(\beta)\Rightarrow \om_1\in\omsort(\alpha)$. Then, in the next step of the algorithm, we apply $\om_2$ on the left to $\alpha'$ and $\beta'$. Since $\alpha'$ is a prefix of $\beta'$, we may use the same reasoning as above to find that $\om_2\in\omsort(\beta)\Rightarrow \om_2\in\omsort(\alpha)$. Continuing in this way we conclude that $\omsort(\alpha)$ is a subword of $\omsort(\beta)$.
\end{proof}

Note that the $\om$-sorting algorithm is defined in terms of {\em left} multiplication by generators; whereas the $\om$-sorting order is comparable to the {\em right} weak order.

\begin{example}
For a more detailed example, we consider the symmetric group $W=\mathfrak{S}_4$ with the generating set of adjacent transpositions,
\begin{equation*}
S=\{s_1=(12),s_2=(23),s_3=(34)\}.
\end{equation*}
Let the sorting word be $\om=(s_1,s_2,s_3,s_2,s_1,s_2)$---another example of the ``cocktail-shaker.'' Since $\om$ is a reduced word for the longest element $w_\circ\in W$, we again have $W_\om=W$. We will use the same notation as in Example \ref{ex:dihedrals} to denote $\om$-sorted words. For example, the string $003210$ denotes the $\om$-sorted word $(s_3,s_2,s_1)$ with index set $\{3,4,5\}\subseteq I(\om)=\{1,2,3,4,5,6\}$. Figure \ref{fig:123212} shows the nested Hasse diagrams of weak order, $\om$-sorting order and Bruhat order.
\begin{figure}
\begin{center}\input{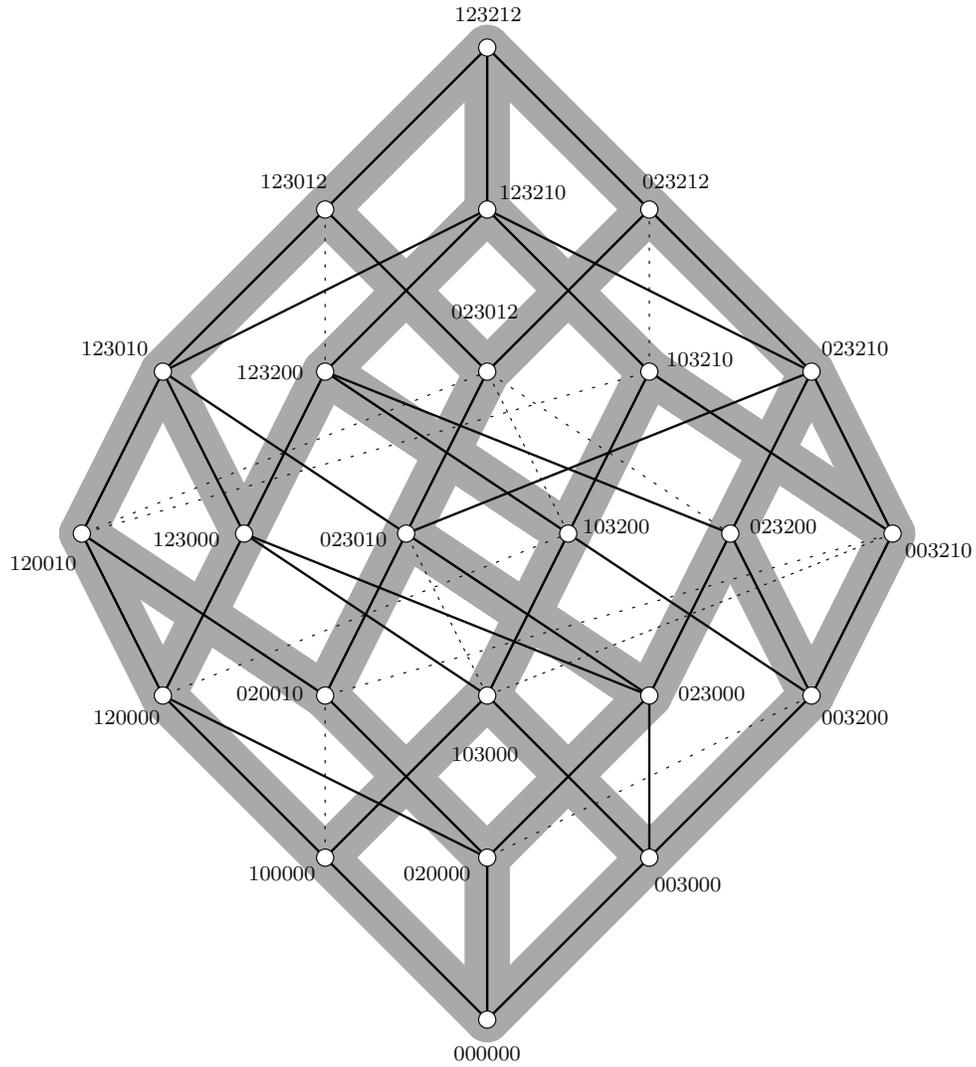}\end{center}
\caption{The weak order, $123212$-sorting order and Bruhat order on $\mathfrak{S}_4$}
\label{fig:123212}
\end{figure}
Weak order is shaded in grey, the solid black lines are the sorting order and the solid plus dotted lines together give Bruhat order.
\end{example}

\subsection{A Supersolvable Antimatroid}
\label{sec:ssantimatroid}
In this section we will show that the collection of $\om$-sorted words forms a supersolvable antimatroid, and hence that the $\om$-sorting order is a supersolvable join-distributive lattice. At the end of the section we will conclude that the sorting orders are parametrized by {\em commutation classes of reduced words}, instead of arbitrary single words.

In addition to the Exchange property we will need the following standard properties of reduced words.

\begin{gluing}
Let $s\in S$ be a generator and consider $\alpha,\beta\in S^*$. If the words $\alpha \beta$, $\alpha s$ and $s\beta$ are all reduced then so is $\alpha s\beta$.
\end{gluing}

\begin{proof}
In this case \cite[Lemma 2.2.10]{bjorner-brenti} implies that $\ell(\langle\alpha\beta\rangle)<\ell(\langle\alpha s\beta\rangle)$. Since $\alpha\beta$ is reduced it follows that $\ell(\langle\alpha s\beta\rangle)=\ell(\langle\alpha\beta\rangle)+1$ and that $\alpha s\beta$ is reduced.
\end{proof}

\begin{lifting}
Given $u\leq_\B w$ and $s\in S$, suppose that $\ell(sw)<\ell(w)$ and $\ell(u)<\ell(su)$---hence, by Exchange we have $sw\leq_\B w$ and $u\leq_\B su$. It follows that $su\leq_\B w$ and $u\leq_\B sw$. See the diagram below.
\begin{center}\input{lifting.pstex_t}\end{center}
\end{lifting}

\begin{proof}
See \cite[Proposition 2.2.7]{bjorner-brenti} for a proof.
\end{proof}

As with the Exchange property, a version of the Lifting property also holds when we multiply on the right by a simple reflection. To see this we replace $u$ and $w$ by $u^{-1}$ and $w^{-1}$, respectively, and note that inversion $w\mapsto w^{-1}$ is an automorphism of Bruhat order.

\begin{theorem}
\label{th:ssantimatroid}
Given an arbitrary Coxeter system $(W,S)$ and a finite sorting word $\om=(\om_1,\ldots,\om_m)\in S^*$, the collection of index sets of $\om$-sorted words,
\begin{equation*}
\F:=\left\{ I(\alpha)\subseteq I(\om): \omsort(\alpha)=\alpha\right\},
\end{equation*}
is a supersolvable antimatroid with respect to the natural order on the ground set $E:=I(\om)=\{1,\ldots,m\}$. (See Definition \ref{def:ssantimatroid}.)
\end{theorem}

\begin{proof}
In the proof we will abbreviate ``lexicographic'' as ``lex.''

To show that $(E,\F)$ is a supersolvable antimatroid, we must verify property \eqref{eq:ssantimatroid}. Consider $\om$-sorted words $\alpha$ and $\beta$ with $\beta\not\subseteq\alpha$, such that $i=\min(I(\beta)\setminus I(\alpha))$. We must show that $\alpha\cup\om_i$ is an $\om$-sorted word---that is, we must show that $\alpha\cup\om_i$ is reduced and that it is the lex-least reduced word for $\langle\alpha\cup\om_i\rangle\in W$ among subwords of $\om$.

In general, $\om_i$ breaks the word $\alpha$ into
\begin{equation*}
\alpha':=\alpha\cap(\om_1,\ldots,\om_{i-1})\quad\text{and}\quad \alpha'':=\alpha\cap(\om_{i+1},\ldots,\om_m),
\end{equation*}
either of which may be empty. Now observe that $\om_i\alpha''$ is a reduced word. Otherwise, by Exchange we obtain a reduced word $\om_i\hat{\alpha}''$ for $\langle\alpha''\rangle$ that is lex-less than $\alpha''$. Consequently, $\alpha'\om_i\hat{\alpha}''$ is a reduced word for $\langle\alpha\rangle$ that is lex-less than $\alpha$, a contradiction. 

Next, we show that $\alpha'\om_i$ is reduced. Set
\begin{equation*}
\beta':=\beta\cap(\om_1,\ldots,\om_{i-1})\quad\text{and}\quad\beta'':=\beta\cap(\om_{i+1},\ldots,\om_m)
\end{equation*}
and note that $\beta'$ is a subword of $\alpha'$ because $i$ is the first place in which $\alpha$ and $\beta$ differ. Furthermore, note that $\beta'\om_i$ is a reduced word since it is a prefix of $\beta$. If $\alpha'\om_i$ is not reduced then the facts $\langle\beta'\rangle\leq_\B\langle\alpha'\rangle$, $\ell(\langle\alpha'\rangle\om_i)<\ell(\langle\alpha'\rangle)$ and $\ell(\langle\beta'\rangle)<\ell(\langle\beta'\rangle\om_i)$, together with the Lifting property, imply that $\langle\beta'\rangle\om_i\leq_\B \langle\alpha'\rangle$. That is, there exists a reduced word for $\langle\beta'\rangle\om_i=\langle\beta'\om_i\rangle$ that is a subword of $\alpha'$. Let $\varphi$ be the lex-least such word. We claim that $\varphi$ is lex-less than $\beta'\om_i$. Indeed, suppose the opposite---let $i'<i$ be the first position in which $\varphi$ and $\beta'\om_i$ differ and suppose that $\om_{i'}\in\beta'\setminus\varphi$. Set
\begin{equation*}
\varphi':=\varphi\cap(\om_1,\ldots,\om_{i'-1})\quad\text{and}\quad\varphi'':=\varphi\cap(\om_{i'+1},\ldots,\om_{i-1}),
\end{equation*}
so that $\varphi=\varphi'\varphi''$. Note that the suffix $\varphi''$ cannot be empty since $\varphi$ and $\beta'\om_i$ have the same number of letters. However if $\varphi''$ is {\em not} empty then by Exchange there exists a subword $\hat{\varphi}''\subseteq\varphi''$ such that $\om_{i'}\hat{\varphi}''$ is a reduced word for $\langle\varphi''\rangle$ and hence $\varphi'\om_{i'}\hat{\varphi}''\subseteq\alpha'$ is a reduced word for $\langle\varphi\rangle$ lex-less than $\varphi$, a contradiction. We conclude that $\varphi$ is lex-less than $\beta'\om_i$ and hence $\varphi\beta''$ is a reduced word for $\langle\beta\rangle$ lex-less than $\beta$. This contradiction proves that the word $\alpha'\om_i$ is reduced.

Since the words $\alpha'\alpha''=\alpha$, $\alpha'\om_i$ and $\om_i\alpha''$ are all reduced, the Gluing property implies that $\alpha'\om_i\alpha''=\alpha\cup\om_i$ is also reduced.

To complete the proof, we must show that $\alpha\cup\om_i$ is lex-least among reduced words for $\langle\alpha\cup\om_i\rangle$. Suppose not. Then by Corollary \ref{cor:sorted} there exists an integer $1\leq j< m$ such that $\om_j$ is not in $\alpha\cup\om_i$ and such that $\gamma:=(\alpha\cup\om_i)\cap(\om_{j+1},\ldots,\om_m)$ has left descent $\om_j$---that is, $\ell(\om_j\langle\gamma\rangle)<\ell(\langle\gamma\rangle)$. But $\delta:=\alpha\cap(\om_{j+1},\ldots,\om_m)$ is a subword of $\gamma$ and we have $\ell(\langle\delta\rangle)<\ell(\om_j\langle\delta\rangle)$. Otherwise, by Exchange we obtain a reduced word $\delta'$ for $\langle\delta\rangle$ that is lex-less than $\delta$, and hence $(\alpha\cap(\om_1,\ldots,\om_{j-1}))\,\delta'$ is a reduced word for $\langle\alpha\rangle$ lex-less than $\alpha$, a contradiction. Since we have $\langle\delta\rangle\leq_\B \langle\gamma\rangle$, $\ell(\om_j\langle\gamma\rangle)<\ell(\langle\gamma\rangle)$ and $\ell(\langle\delta\rangle)<\ell(\om_i\langle\delta\rangle)$, the Lifting property tells us that $\om_j\langle\delta\rangle\leq_\B \langle\gamma\rangle$. That is, there exists a subword $\gamma'\subseteq\gamma$ that is a reduced word for $\om_j\langle\delta\rangle$. In this case, $(\alpha\cap(\om_1,\ldots,\om_{j-1}))\,\om_j\gamma'$ is a reduced word for $\langle\alpha\rangle$ that is lex-less than $\alpha$. This final contradiction proves that $\alpha\cup\om_i$ is lex-least among reduced words for $\langle\alpha\cup\om_i\rangle$, as we wished to show.
\end{proof}

Following this, Theorem \ref{th:ssjd} implies our main result. Recall that a join-distributive lattice is graded by the cardinality of feasible sets in the corresponding antimatroid.

\begin{corollary}
The $\om$-sorting order on $W_\om\subseteq W$ is a supersolvable join-distributive lattice, graded by the usual Coxeter length function $\ell:W\to\Z$.
\end{corollary}

We also conclude that the sorting word might as well be reduced.

\begin{corollary}
\label{cor:reduced}
We lose nothing if we consider only reduced sorting words---that is, given a sorting word $\om\in S^*$ there exists a reduced word $\om'\subseteq\om$ for which $\om'$-sorting is the same as $\om$-sorting.
\end{corollary}

\begin{proof}
Given an arbitrary word $\om\in S^*$, let $\om'\subseteq\om$ denote the union of all $\om$-sorted subwords of $\om$. Since the collection of feasible sets of an antimatroid is closed under taking unions (Property \eqref{eq:unions}), we conclude that $\om'$ is an $\om$-sorted word and hence it is reduced. Since the elements of $\om\setminus\om'$ will play no role in the antimatroid, nor in the sorting order, we may replace $\om$ by the reduced word $\om'$ without loss.
\end{proof}

Furthermore, if $\om'\subseteq\om$ denotes the union of all $\om$-sorted words, we note that the collection of elements $W_\om\subseteq W$ that occur as subwords of $\om$ coincides with the lower interval $[1,\langle\om'\rangle]_\B$ in Bruhat order. Indeed, any element $w\in W_\om$ has a reduced $\om$-sorted word, which by definition is a subword of the reduced word $\om'$. Thus we may think of the sorting orders as partial orders on a lower interval $[1,w]_\B$ in Bruhat order, parametrized by reduced word for $w$.

Finally we note that the $\om$-sorting order depends only on the ``commutation class'' of the word $\om$. We say that two words are in the same {\sf commutation class} if one can be obtained from the other by repeatedly exchanging adjacent commuting generators.
 
\begin{lemma}
Let $\om,\zeta\in S^*$ be two words that differ by the exchange of commuting generators in positions $i$ and $i+1$. Then we have
\begin{equation*}
u\leq_\om w\quad\Longleftrightarrow\quad u\leq_\zeta w
\end{equation*}
for all $u,w\in W_\om=W_\zeta\subseteq W$.
\end{lemma}

\begin{proof}
Consider $\om=(\om_1,\ldots,\om_m)$ and $\zeta=(\zeta_1,\ldots,\zeta_m)$, where $\om_i=\zeta_{i+1}$, $\om_{i+1}=\zeta_i$ and $\om_i\om_{i+1}=\om_{i+1}\om_i$. Note that the transposition $(i,i+1)$ acts on subsets of $\{1,\ldots,m\}$ by switching the indices $i$ and $i+1$. If $\alpha$ is a subword of $\om$ with index set $I(\alpha)\subseteq I(\om)=\{1,\ldots,m\}$, let $\alpha'$ denote the subword of $\zeta$ with index set $(i,i+1)\cdot I(\alpha)$. We claim that the involution $\alpha\mapsto \alpha'$ (which satisfies $\langle\alpha\rangle=\langle\alpha'\rangle$) is a bijection between $\om$-sorted and $\zeta$-sorted words, from which the result follows.

Indeed, suppose $\alpha$ is $\om$-sorted. We wish to show that $\alpha'$ is $\zeta$-sorted. (The proof that $\alpha'\mapsto\alpha$ preserves sortedness will be the same.) Since we have merely exchanged commuting generators, $\alpha'$ is reduced. We must show that $\alpha'$ is the lex-least reduced word for $\langle\alpha'\rangle$ among subwords of $\zeta$. Suppose not. Then by Corollary \ref{cor:sorted} there exists $1\leq j<m$ such that $\zeta_j\not\in\alpha'$ and $\zeta_j$ is a left descent for $\alpha'\cap(\zeta_{j+1},\ldots,\zeta_m)$. If both or neither of $\zeta_i,\zeta_{i+1}$ occur in $\alpha'$, or if $j\not\in\{i,i+1\}$, we find that $\om_j\not\in\alpha$ and $\om_j=\zeta_j$ is a left descent of
\begin{equation*}
\langle\alpha\cap(\om_{j+1},\ldots,\om_m)\rangle=\langle\alpha'\cap(\zeta_{j+1},\ldots,\zeta_m)\rangle,
\end{equation*}
which implies that $\alpha$ is not lex-least, a contradiction. Otherwise, exactly one of $\zeta_i,\zeta_{i+1}$ occurs in $\alpha$---without loss of generality, say $\zeta_i\in\alpha'$---and we have $j=i+1$. In this case, $\om_i\not\in\alpha$ and $\om_i=\zeta_{i+1}$ is a left descent of
\begin{equation*}
\langle\alpha\cap(\om_{i+2},\ldots,\om_m)\rangle=\langle\alpha'\cap(\zeta_{i+2},\ldots,\zeta_m)\rangle.
\end{equation*}
Since $\om_i$ and $\om_{i+1}$ commute, $\om_i$ is also a left descent of
\begin{equation*}
\om_{i+1}\langle\alpha\cap(\om_{i+2},\ldots,\om_m)\rangle=\langle\alpha\cap(\om_{i+1},\ldots,\om_m)\rangle,
\end{equation*}
which implies that $\alpha$ is not lex-least, again a contradiction.
\end{proof}

\subsection{Remark---A Maximal Lattice}
Recall that Bruhat order is obtained from weak order by adding the extra cover relations of the form $u\prec w$ where $w=ut$ for some non-simple reflection $t\in T\setminus S$ such that $\ell(w)=\ell(u)+1$. In general, let $(P,\leq)$ be a graded poset. We say that another graded poset $(P,\leq')$ is a {\sf graded extension} of $(P,\leq)$ if $x\leq y$ implies $x\leq' y$ for all $x,y\in P$ and if the rank function is the same for both posets. In this case we also say that $(P,\leq)$ is a {\sf graded contraction} of $(P,\leq')$ We are interested in the collection of graded extensions between the weak and Bruhat orders and the role that the sorting orders play among these.

Our main observation\footnote{Thanks to Hugh Thomas for suggesting that our original observation for one cover could be generalized to any number of covers.} is the following.

\begin{definition}
Let $(P,\leq)$ be a finite graded lattice with rank function $\rk:P\to\Z$. We say that $P$ is {\sf maximal} if the addition of any finite collection of cover relations of the form $x\prec y$ with $\rk(y)=\rk(x)+1$ yields a nonlattice.
\end{definition}

It turns out that all supersolvable join-distributive lattices---hence, in particular, all distributive lattices---are maximal. The converse of this theorem is not true.

\begin{theorem}
Let $(P,\leq)$ be a supersolvable join-distributive lattice. Then $P$ is maximal.
\end{theorem}

\begin{proof}
By Theorem \ref{th:ssjd}, $P$ is the lattice of feasible sets of a supersolvable antimatroid $(E,\F)$. Consider a collection $\{ (A_i,B_i)\in\F^2\}$ of pairs of feasible sets such that $\abs{B_i}=\abs{A_i}+1$ and $A_i\not\subseteq B_i$ for all $i$, and choose $j$ such that the cardinality of $A_j$---and hence $B_j$---is minimal.

Let $P'$ denote the poset obtained from $P$ by adding the cover relations $A_i\prec B_i$ for all $i$. We claim that $P'$ is not a lattice. Indeed, let $C_j:=A_j\cup\{x_j\}$ where $x_j=\min(B_j\setminus A_j)$. By Property \eqref{eq:ssantimatroid} we have $C_j\in\F$. Note that the meet of $C_j$ and $B_j$ in $P$ is given by the union of all feasible sets contained in $C_j\cap B_j$:
\begin{equation*}
C_j\wedge B_j = \bigcup_{\substack{ X\in\F \\ X\subseteq C_j\cap B_j }} X.
\end{equation*}
Since the prefix of every feasible set is feasible, so is the prefix of $B_j$ ending in $x_j$. By construction this prefix is also contained in $C_j$, hence $x_j\in C_j\wedge B_j$. Thus $A_j$ and $C_j\wedge B_j$ are two lower bounds for $C_j$ and $B_j$ in $P'$ that are incomparable in $P$. Since the cardinality of $A_j$ was chosen to be minimal, $A_j$ and $C_j\wedge B_j$ are also incomparable in $P'$. We conclude that $P'$ is not a lattice.
\end{proof}

In particular, the $\om$-sorting order is a maximal lattice with respect to the addition of any collection of Bruhat cover relations of the form $u\prec ut$.

\begin{figure}
\begin{center}\input{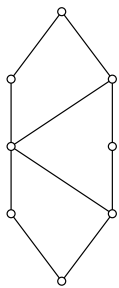}\end{center}
\caption{A maximal lattice that is not a sorting order}
\label{fig:maxlattice}
\end{figure}

To extend this result, one might try to classify all maximal lattices among the graded extensions between weak and Bruhat order---or more generally among all graded contractions of Bruhat order. Note that the sorting orders do not provide the complete solution to this problem, since the lattice in Figure \ref{fig:maxlattice} is maximal between the weak and Bruhat orders on the dihedral group $I_2(4)$ (a.k.a. $B_2$), but it is not a sorting order in our sense.

Our study of the general properties of $\om$-sorting order is now complete. In the remaining sections we consider two special cases---those of {\em infinite} sorting words and {\em cyclic} sorting words.

\section{Infinite Sorting Orders}
\label{sec:infinite}

Thus far we have dealt exclusively with finite sorting words because some of the structures in Section \ref{sec:convexity} make sense only in the finite case. However, most of our results in this paper still hold when $\om$ is a semi-infinite word,
\begin{equation*}
\om=(\om_1,\om_2,\om_3,\ldots).
\end{equation*}
For example, the characterization of $\om$-sorted words via the sorting algorithm and the definition of $\om$-sorting order remain the same.

If there are only finitely many group elements that occur as subwords of $\om$ (which happens, for instance, when $W$ is a finite group), then we may restrict $\om$ to the (finite) union of all $\om$-sorted words and everything in the paper goes through as before. Hence, let us assume that $W$ is infinite and that infinitely many group elements occur as subwords of the sorting word $\om$.

In this case, Theorem \ref{th:weakBruhat} still holds---the $\om$-sorting order extends the weak order and is extended by Bruhat order. Theorem \ref{th:ssantimatroid} also goes through as before---the collection of $\om$-sorted words forms a supersolvable antimatroid $(E,\F)$ on an infinite ordered ground set $E$. Although Korte, Lov\'asz and Schrader \cite{KLS} did not consider antimatroids with infinite ground sets, the conclusions of Lemma \ref{lem:antimatroid} still hold in this case.\footnote{Infinite antimatroids and convex geometries are considered in \cite{AGT}.} Thus, we have the following.

\begin{definition} 
Let $(P,\leq)$ be a locally-finite lattice (all intervals are finite). We say that $P$ is {\sf join-distributive} if every atomic interval in $P$ is boolean.\footnote{This terminology is not entirely standard; see \cite{AGT}.} We say that $P$ is {\sf supersolvable} if every interval in $P$ is supersolvable in the usual sense.
\end{definition}

\begin{theorem}
Let $(W,S)$ be an infinite Coxeter system and let $\om$ be a word in which infinitely many group elements occur as subwords. Then the $\om$-sorting order $(P_\om,\leq_\om)$ is an infinite lattice that is supersolvable and join-distributive.
\end{theorem}

\begin{proof}
The poset $P_\om$ contains joins because the collection of feasible sets $\F$ is closed under taking unions (Property \eqref{eq:unions}). Then since $P_\om$ has a minimum element (the identity $1_W\in W$), it also contains meets. The other properties follow as before.
\end{proof}

This conclusion is remarkable because the weak order on an infinite Coxeter group is {\em not} a lattice---while is possesses meets, it does not possess joins---and Bruhat order is not a lattice even in the finite case. It is thus interesting to have a new source of lattice structures on the elements of an infinite Coxeter group. Indeed, we know of no other source.

It also remains true that the sorting order only depends on the commutation class of the sorting word and we may still assume that the sorting word is reduced in the following sense.

\begin{definition}
\label{def:reduced}
We say that an infinite word $\om$ in the generators $S$ is {\sf reduced} if every prefix of $\om$ is reduced in the usual sense.
\end{definition}

Indeed, we may restrict $\om$ to the union $\om'\subseteq\om$ of all $\om$-sorted words. Any prefix $\alpha$ of $\om'$ is then contained in the union of {\em finitely} many $\om$-sorted words, which is reduced. Since $\alpha$ is the prefix of a reduced word, it is reduced.

\section{Cyclic Sorting}
\label{sec:cyclic}

Finally, we discuss an important special case of $\om$-sorting which has been the motivation for our work. Let $(W,S)$ be an arbitrary Coxeter system with generators $S=\{s_1,s_2,\ldots,s_n\}$. Any word of the form
\begin{equation*}
c=(s_{\sigma(1)},s_{\sigma(2)},\ldots,s_{\sigma(n)}),
\end{equation*}
where $\sigma$ is a permutation of $\{1,\dots,n\}$, is called a {\sf Coxeter word}; the corresponding group element $\langle c\rangle\in W$ is a {\sf Coxeter element}. We say that a {\sf cyclic word} is any semi-infinite word of the form
\begin{equation*}
c^\infty:=ccc\ldots,
\end{equation*}
where $c$ is a Coxeter word.\footnote{In the case that every irreducible component of $W$ is infinite, Speyer \cite{speyer} has recently shown that cyclic words are reduced in the sense of Definition \ref{def:reduced}.} In this case Reading was the first to consider the $c^\infty$-sorting algorithm---which he called ``$c$-sorting''---for elements of the group $W$ (see \cite{reading:csort1,reading:csort2}). However he did not consider the structure of the collection of sorted words nor the corresponding partial order.

His main interest was the collection of so called ``$c$-sortable'' elements. Because of the cyclic nature of $c^\infty$, each $c^\infty$-sorted word $\alpha$ naturally splits into a sequence of subwords $\alpha_{(1)}\alpha_{(2)}\alpha_{(3)}\ldots$, where $\alpha_{(i)}$ is the intersection of $\alpha$ with the index set $\{(i-1)n+1,(i-1)n+2,\ldots,in\}$. Given a word $\beta\in S^*$, let $\tilde{\beta}\subseteq S$ denote its underlying set of letters.

\begin{definition}
Given a group element $w\in  W$, let $\alpha$ denote its $c^\infty$-sorted word. We say that the element $w$ is {\sf $c$-sortable} if we have a descending chain
\begin{equation*}
\tilde{\alpha}_{(1)}\supseteq \tilde{\alpha}_{(2)}\supseteq \tilde{\alpha}_{(3)}\supseteq\cdots
\end{equation*}
of subsets of $S$.
\end{definition}

Reading introduced this notion because the $c$-sortable elements are precisely the elements of his $c$-Cambrian lattice \cite[Theorem 1.1]{reading:csort2}. Probably their most remarkable property is contained in the following theorem.

\begin{theorem}{\bf \cite[Theorem 9.1]{reading:csort1}} Let $(W,S)$ be a {\em finite} Coxeter system with $\left| S\right|=n$. Given any Coxeter word $c$, the number of $c$-sortable elements in $W$ is equal to
\begin{equation*}
{\rm Cat}(W)=\prod_{i=1}^n \frac{h+d_i}{d_i},
\end{equation*}
where $h$ is the Coxeter number (the order of a Coxeter element) and $\{d_1,d_2,\ldots,d_n\}$ is the multiset of degrees of fundamental invariants for $(W,S)$.
\end{theorem}

This ``generalized Catalan number'' ${\rm Cat}(W)$ has played a central role in much recent work (for references see \cite{armstrong}), and the $c$-sortable elements provided Reading with a bridge between two important classes of ``Catalan objects'': the clusters and the noncrossing partitions.

It is natural to ask what special properties the sorting order has in the case of a cyclic sorting word $c^\infty$. The following theorem provides a partial answer.

\begin{theorem}
Let $c^\infty$ be a cyclic sorting word for an arbitrary Coxeter system $(W,S)$. The $c^\infty$-sorting order restricted to $c$-sortable elements is a join-distributive lattice.
\end{theorem}

\begin{proof}
The proof follows from two observations. First note that the $c$-sortable elements form a join-sublattice of the full $c^\infty$-sorting order. This is because the defining property of $c$-sortability is preserved under taking unions of words.

Second, let $u$ and $w$ be $c$-sortable group elements with $c^\infty$-sorted words $\alpha$ and $\beta$, respectively, and suppose that $u\prec w$ is a cover in the sorting order on $c$-sortable elements. We claim that this is also a cover in the full $c^\infty$-sorting order. Suppose not, so that $\ell(w)>\ell(u)+1$. If $x=\min(\beta\setminus\alpha)$, then $\langle\alpha\cup x\rangle$ is another $c$-sortable group element strictly between $u$ and $w$, contradicting the fact that $u\prec w$ is a cover.

Finally, it is easy to see that a join-sublattice of a join-distributive lattice that preserves covers is also a join-distributive lattice. 
\end{proof}

In a forthcoming paper we will show that the $c^\infty$-sorting and Bruhat orders coincide on $c$-sortable elements, and moreover that this order is supersolvable. For now we present an example.

\begin{figure}
\begin{center}\input{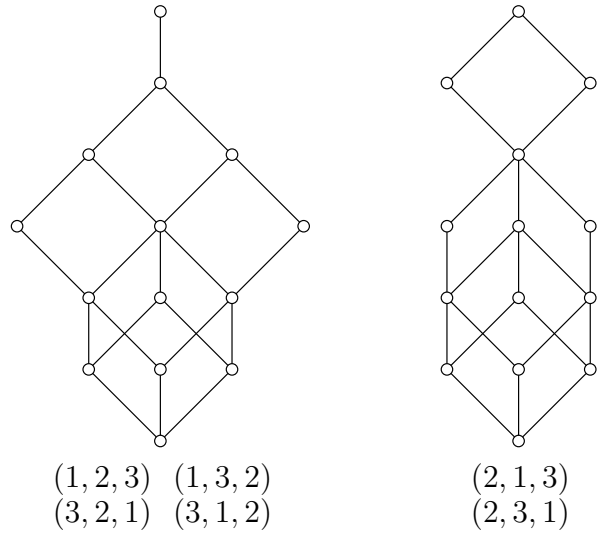}\end{center}
\caption{The two possible lattices of sortable elements in type $A_3$}
\label{fig:sortables}
\end{figure}

\begin{example}
Consider the Coxeter system $(W,S)$ of type $A_3$ with Coxeter diagram:
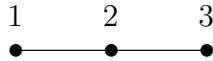
\begin{center}\input{a3diagram.pstex_t}\end{center}
Among the six possible Coxeter words $c$, there are just two possibilities for the isomorphism type of the lattice of $c$-sortable elements under $c^\infty$-sorting order. These are displayed in Figure \ref{fig:sortables} with the corresponding Coxeter words. Notice that one of these is the well-known lattice of order ideals of the root poset. This phenomenon, unfortunately, does not persist for all types.
\end{example}

\section*{Acknowledgments}
The author gratefully acknowledges helpful conversations with Peter McNamara, Nathan Reading, Vic Reiner, David Speyer and Hugh Thomas, as well as the helpful comments of the anonymous referees.

\bibliographystyle{amsalpha}

\begin{thebibliography}{A}

\bibitem{AGT}
K.~Adaricheva, V.~Gorbunov and V.~Tumanov, \emph{Join-semidistributive lattices and convex geometries}, Advances in Math. {\bf 173} (2003), 1--49.

\bibitem{armstrong}
D.~Armstrong, \emph{Generalized noncrossing partitions and combinatorics of Coxeter groups}, \texttt{arxiv:math.CO/0611106}, to appear in Mem. Amer. Math. Soc.

\bibitem{bjorner-brenti}
A.~Bj\"orner and F.~Brenti, \emph{Combinatorics of Coxeter groups}, Springer (2005).

\bibitem{dilworth}
R.~Dilworth, \emph{Lattices with unique irreducible decompositions}, Ann. of Math. {\bf 41} (1940), 771--777.

\bibitem{edelman:meet}
P.~Edelman, \emph{Meet-distributive lattices and the anti-exchange closure}, Algebra Universalis {\bf 10} (1980), 290--299.

\bibitem{edelman-jamison}
P.~Edelman and R.~Jamison, {\em The theory of convex geometries}, Geometriae Dedicata {\bf 19} (1985), 247--270.

\bibitem{gratzer}
G.~Gr\"atzter, \emph{General lattice theory}, Academic Press (1978).

\bibitem{HR}
M.~Hawrylycz and V.~Reiner, \emph{The lattice of closure relations on a poset}, Algebra Universalis {\bf 30} (1993), 301--310.

\bibitem{knuth:taocp}
D.~Knuth, \emph{The art of computer programming, Volume 1: Fundamental algorithms}, Addison-Wesley (1973).

\bibitem{knuth:aah}
D.~Knuth, \emph{Axioms and hulls}, Lecture Notes in Computer Science, no. 606, Springer (1992).

\bibitem{knutson-miller}
A.~Knutson and E.~Miller, \emph{Subword complexes in Coxeter groups}, Advances in Math. {\bf 184} (2004), 161--176.

\bibitem{KLS}
B.~Korte, L.~Lov\'asz and R.~Schrader, {\em Greedoids}, Algorithms and Combinatorics {\bf 4}, Springer (1991).

\bibitem{mcnamara}
P.~McNamara, \emph{\EL-labellings, supersolvability and $0$-Hecke algebra actions on posets}, J. Combin. Theory Ser. A {\bf 101} (2003), 69--89.

\bibitem{reading:csort1}
N.~Reading, \emph{Clusters, Coxeter-sortable elements and noncrossing partitions}, Trans. Amer. Math. Soc. {\bf 359} (2007), 5931--5958.

\bibitem{reading:csort2}
N.~Reading, \emph{Sortable elements and Cambrian lattices}, Algebra Universalis {\bf 56} (2007), 411--437.

\bibitem{speyer}
D.~Speyer, \emph{Powers of Coxeter elements in infinite groups are reduced}, \texttt{arXiv:0710.3188}

\bibitem{stanley:ec1}
R.~Stanley, \emph{Enumerative Combinatorics}, vol. 1, Cambridge University Press (1997).

\bibitem{stanley:ss}
R.~Stanley, \emph{Supersolvable lattices}, Algebra Universalis {\bf 2} (1972), 197--217.

\end{thebibliography}

\end{document}